\documentclass[a4paper,12pt]{article}

\usepackage{amsmath}
\usepackage{amsthm}
\usepackage{amssymb}
\usepackage{latexsym}
\usepackage{graphicx, color}
\usepackage{stmaryrd}
\usepackage{mathrsfs}
\usepackage{enumerate}
\usepackage{cases}

\font\sy=rsfs10     scaled 1200

\newcommand{\glalign}[2]{\lower.6ex\vbox{
\baselineskip\lineskip\ialign{$#1\hfil##\hfil$\crcr#2\crcr=\crcr}}}

\newcommand{\scr}[1]{{\mbox{\sy #1}\,}}

\newcommand{\del}{\partial}

\renewcommand{\div}{\mbox{\rm div}\,}

\newcommand{\supp}{\mbox{\rm supp}\,}

\newcommand{\trans}{{}^\top}

\def\eqn#1$$#2$${\begin{equation}\label#1#2\end{equation}}

\numberwithin{equation}{section}

\newtheorem{defi}{Definition}[section]
\newtheorem{thm}[defi]{Theorem}

\newtheorem{prop}[defi]{Proposition}
\newtheorem{lem}[defi]{Lemma}
\newtheorem{rem}[defi]{Remark}
\def\eqn#1$$#2$${\begin{equation}\label#1#2\end{equation}}

\numberwithin{equation}{section}
\numberwithin{equation}{section}
\allowdisplaybreaks[4]

\begin{document}

%------------------------------------------------------------
%     Section 1.   
%------------------------------------------------------------

\title{
\bf \large 
Global existence and time decay estimate of solutions to the compressible Navier-Stokes-Korteweg system under critical condition}
\author{
Takayuki KOBAYASHI and Kazuyuki TSUDA\\
 Osaka University, \\
 1-3, Machikaneyamacho, Toyonakashi, 560-8531, JAPAN \\
e-mail: jtsuda@sigmath.es.osaka-u.ac.jp}
\date{}
\maketitle
%%%%%%%%%%%%%%%%%%%%%%
\begin{abstract} 
Global existence of solutions to the compressible Navier-Stokes-Korteweg system around a constant state is studied. This system describes liquid-vapor two phase flow with phase 
transition as diffuse interface model. In previous works they assume that the pressure is a monotone function for change of density similarly to the usual compressible Navier-Stokes system. On the other hand, 
due to phase transition the pressure is accurately non-monotone function and the linearized system loses symmetry in a critical case such that the derivative of pressure is 0 at the given constant state. 
It is shown that  in the critical case for small data whose momentum has derivative form there exist global $L^2$ solutions and the parabolic type decay rate of the solutions is obtained.
The proof is based on decomposition method for solutions to a low frequency part and a high frequency part.  
\end{abstract}



%------------------------------------------------------------
%     Section 1.   
%------------------------------------------------------------

\section{Introduction}
We study global existence of solutions to the following compressible Navier-Stokes-Korteweg system 
in $\mathbb{R}^n$ $(n \geq 3)$:
\begin{numcases}
{}
\partial_{t}\rho +\div m=0,\nonumber\\
\partial_{t}m +\div \Big(\frac{m \otimes m}{\rho}\Big)+\nabla P(\rho)=\div \Big({\mathcal S}(\frac{m}{\rho})+{\mathcal K}(\rho)\Big),\label{CNSK}\\
\rho(x,0) =\rho_0, \ \ m(x,0)=m_0. \nonumber
\end{numcases}
\noindent Here $\rho=\rho(x,t)$ and $m=(m_{1}(x,t),\cdots,m_{n}(x,t))$ denote the unknown density and momentum respectively, at time $t\in \mathbb{R}_{+}$ and 
position $x\in\mathbb{R}^n$; $\rho_0=\rho_0(x)$ and $m_0=m_0(x)$ denote given initial data; 
${\mathcal S}$ and ${\mathcal K}$ denote the viscous stress tensor and the Korteweg stress tensor that are given by 
\begin{eqnarray}
\left\{
\begin{array}{ll}
{\mathcal S}(\frac{m}{\rho}) =\Big(\mu' \div \frac{m}{\rho} \Big)\delta_{i,j} +2\mu d_{ij}\Big(\frac{m}{\rho}\Big),\\
{\mathcal K}(\rho) =\frac{\kappa}{2}(\Delta \rho^2 -|\nabla \rho|^2)\delta_{i,j}- \kappa \frac{\del \rho}{\del x_i}\frac{\del \rho}{\del x_j},
\end{array}
\right.\label{Korteweg tensor}
\end{eqnarray}
where $d_{ij}\Big(\frac{M}{\rho}\Big)=\frac{1}{2}\left(\frac{\del }{\del x_i}\Big(\frac{M}{\rho}\Big)_j+\frac{\del}{\del x_j}\Big(\frac{M}{\rho}\Big)_i\right)$; $\mu$ and $\mu'$ are the viscosity coefficients that are assumed to be constants satisfying 
$$
\mu>0, \ \ \ 
\frac{2}{n}\mu+\mu'\geq 0. 
$$   
$\kappa$ denotes the capillary constant that is assumed to be a positive constant. Note that if $\kappa=0$ in the Korteweg tensor, the usual compressible Navier-Stokes equation (the abbreviation is used by ``CNS'' below) is obtained; 
$P=P(\rho)$ is the pressure that is assumed to be a smooth function of $\rho$. Here we assume that $P$ satisfies   
\begin{eqnarray}
P'(\rho_*)=0, \label{nonsymmetric}
\end{eqnarray}
where $\rho_*$ is a given positive constant and $(\rho_*,0)$ denotes a given constant state. We consider solutions to \eqref{CNSK} around the constant state.

\eqref{CNSK} governs motion of two phase flow
between liquid and vapor with phase transition in a compressible fluid. 
To describe the phase transition, this system use diffuse interface. The phase boundary is regarded as a narrow transition layer and the fluid state is described by a phase parameter, change of the density in this system. Therefore it is enough to analysis one set of equations in a single spatial domain. Furthermore difficulty of topological change of the interface does not occur in difference from the classical sharp interface model. 
Van der Waals \cite{Van der Waals} suggests diffuse interface model which occurs from a steep gradient of the density for the liquid-vapor type two phase flow. Based on his idea, Korteweg \cite{Korteweg}  modifies the stress tensor of the usual Navier-Stokes equation.  The modified stress tensor includes $\nabla\rho \otimes \nabla \rho$ similarly to \eqref{Korteweg tensor}. 
Dunn and Serrin \cite{Dunn and Serrin} generalize the Korteweg's work and derive the system $(\ref{CNSK})$ with (\ref{Korteweg tensor}) rigorously. 
Heida and M\'{a}lek \cite{Heida and malek} 
also derive \eqref{CNSK} by  the entropy production method which does not require to introduce any new or non-standard concepts such as multipolarity or interstitial working which are used in \cite{Dunn and Serrin}.

Our aim is to show global existence of solutions to \eqref{CNSK} and study convergence rate of the solution to the given constant state under the condition \eqref{nonsymmetric}. 
Concerning global existence of solutions to \eqref{CNSK} on $\mathbb{R}^n$, as far as we investigated, all study assume that $P'(\rho_*)>0$, which is the same condition as that of CNS \cite{Matsumura-Nishida}.  
Concretely, Hattori and Li \cite{Hattori-Li-1,Hattori-Li-2} obtain the global existence of $H^{N+1} \times H^{N}$ solutions  with a small initial data  
$u_0 \in H^{N+1}\times H^N$, where $H^N$ denotes the usual $L^2$ Sobolev space and $N$ is an integer satisfying that $N\geq [n/2]+2$ and $[n/2]$ denotes the integer part of $n/2$. Danchin and Desjardins \cite{Danchin} show the global existence of solutions with small initial data $u_0 \in (B^{\frac{n}{2}}_{2,1}\cap B^{\frac{n}{2}-1}_{2,1})\times B^{\frac{n}{2}-1}_{2,1}$, where $B^{\frac{n}{2}}_{2,1}$ denotes the usual homogeneous Besov space. 
Recently, Tan and R. Zhang, X. Zhang and Tan and Tan, Wang and Xu \cite{Tan-Zhang,Zhang-Tan, Tan-Wang-Xu} show the global existence for small initial data in some Sobolev spaces which have lower regularity that that of \cite{Hattori-Li-1,Hattori-Li-2} in three dimensional case. 
In addition, Wang and Tan \cite{Wang-Tan} study convergence rates of $L^p$ $(2\leq p)$ norms of the solutions. 
They show that if initial data satisfy $\|(\rho_0, v_0)\|_{(H^{s+1}\times H^s)\cap L^1}<<1$ $(s\geq 3)$, where $v$ denotes the velocity field $v=\frac{m}{\rho}$, it holds hat for $t>0$ 
\begin{eqnarray*}
&&\|(\rho(t)-\rho_*, v(t))\|_{L^p}\leq C(1+t)^{-\frac{3}{2}(1-\frac{1}{p})}  \ \ (2\leq p \leq 6),\\[1ex]
&&\|\nabla (\rho(t)-\rho_*, v(t))\|_{L^2}\leq C(1+t)^{-\frac{5}{4}}. 
\end{eqnarray*}

However, as shown in J. Daube \cite{de-bu} the pressure is non-monotone function due to the phase transitions.  
Indeed, the pressure $P$ is given by the Van der Waals equation of state
$$
P(\rho) = \rho^2 \varphi'(\rho)
$$
for a given smooth specific Helmholtz energy $\varphi=\varphi(\rho)$. In order to model phase transitions, it is assumed that the Helmholtz free energy $\tilde{W}(\rho)=\rho\varphi(\rho)$ has a double-well sharp. (Figure.1) Hence, as shown in \cite{de-bu}, this together with the relation between $P$ and $\tilde{W}(\rho)$; 
$$
P(\rho)=\rho\tilde{W}'(\rho)-\tilde{W}(\rho)
$$
show that the pressure is a non-monotone function of the density. (Figure.2).

\begin{figure}[h]
\begin{center} 
%\centerline{\includegraphics{fx1}\hspace*{5mm}\includegraphics{fx1}}
\begin{minipage}{0.5\hsize}
\centerline{\includegraphics[height=4.0cm, width=4.0cm]{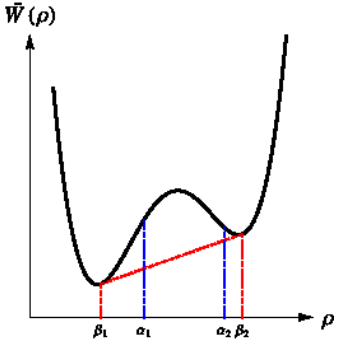}}
\caption{Typical shape of $\tilde{W}=\tilde{W}(\rho)$(\cite{de-bu})
}
\label{w} \ \ \ \ 
\centerline{\includegraphics[height=6.0cm, width=8.0cm]{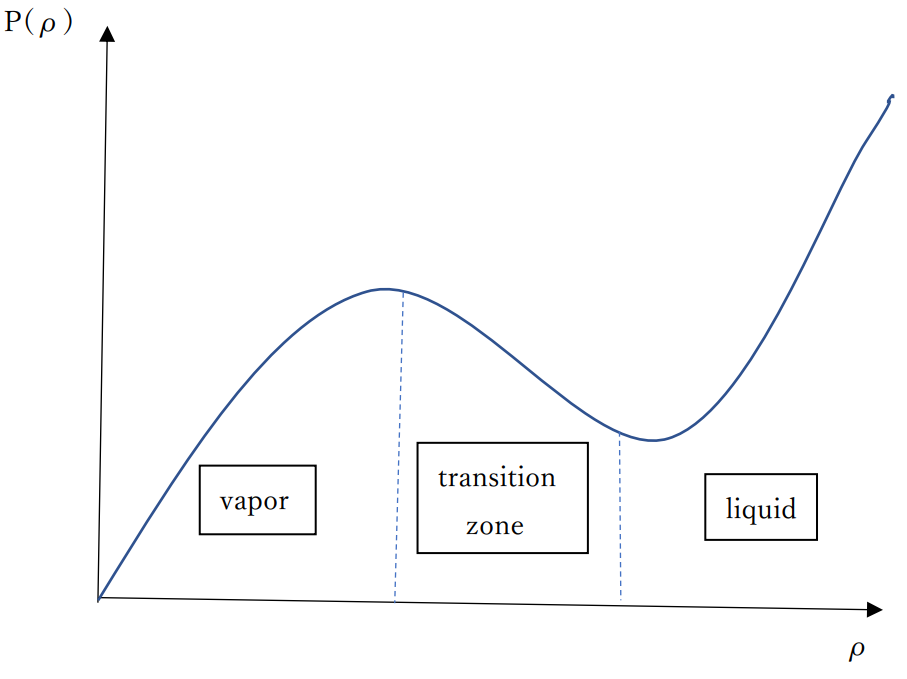}}
\caption{Typical shape of the pressure $P=P(\rho)$}\vspace*{-6pt}
\label{pressure}
\end{minipage}
\end{center}
\end{figure}

Consequently, we should consider the case not only $P'(\rho_*)>0$ but also $P'(\rho_*) \leq 0$ in \eqref{CNSK}. 
When $P'(\rho_*)<0$, Figure. 2 shows that the fluid state is in the phase transition which is mixture state between liquid and vapor. Hence we can not expect that the constant state is stable and we have global existence of solutions around the constant state. On the other hand if $P'(\rho_*)>0$ it is already mentioned above that  the many previous results show that the constant state is stable.  
These motivate us to study that in the critical case $P'(\rho_*) =0$, whether we have the global existence of solution around the constant state or not. 
Mathematical difficulty is that when $P'(\rho_*) =0$, the linear system loses symmetry for a linear derivative operator of spacial variables even if we assume that $\kappa=0$. It is well known that the symmetry has an basic role for stability condition and $L^2$ decay estimate of solutions in general hyperbolic conservation system with relaxation terms which includes CNS as in Kawashima, Shizuta and Umeda and Kawashima  and Shizuta \cite{Kawashima1, Kawashima2}. 
Furthermore, due to $P'(\rho_*)=0$ the momentum part of the fundamental solutions to linear system in a low frequency  has worse order terms in the Fourier space as shown in \eqref{solution-linear} below than that of Kobayashi and Shibata \cite{Kobayashi-Shibata} for linearized CNS. This fact prevents us from getting the parabolic type time decay estimate of solutions. 

We shall show that for \eqref{CNSK} there exist global $L^2$ solutions for small data $u_0=\trans(\rho_0,m_0)$ with  
a regularity assumption such that $m_0$ has the derivative form $m_0= \del_x \tilde{m}_0$. Furthermore, the solutions converge to the constant state with the parabolic type decay rate;
$$
\|\nabla^k(\rho(t)-\rho_*, m(t))\|_{L^2}\leq C(1+t)^{-\frac{n}{4}-\frac{k}{2}},
$$
for $k=0,1$. This rate coincides with those of \cite{Wang-Tan} for the case $P'(\rho_*)>0$ and \cite{Matsumura-Nishida} for CNS in the three dimensional space.  

To show the global existence theorem, we introduce decomposition of solution to a low frequency part and a high frequency part as in Okita \cite{Okita} and Tsuda \cite{Tsuda-CNSK}.  
Concerning linear estimate of the low frequency part, we use a similar method to that of \cite{Kobayashi-Shibata} for CNS and Shibata \cite{Shibata-v} for the linear viscoelastic system.  
By virtue of combining explicit forms of the fundamental solutions with the regularity assumption, we overcome worse order terms which appear in the density part than that of \cite{Kobayashi-Shibata} and we can get the estimate with the same decay order as that of solutions to the heat equation. Note that by using the conservation form, the nonlinearity satisfies the regularity assumption and we can also estimate the nonlinear problem. 

On the other hand, as for the high frequency part, we use $L^2$ energy method in the Fourier space. 
Since the linear system loses the symmetry in the conservation of momentum, to obtain closed estimate for $L^2$  norm of the density is a key point. To get the closed estimate we combine Hattori and Li \cite{Hattori-Li-1,Hattori-Li-2} type $L^2$ energy method and the Poincar\'{e} type estimate which holds in the high frequency part. Then
we can derive the linear estimate for both the density and momentum in the usual $L^2$ Sobolev spaces. 
Concerning nonlinear estimates, note that $\rho$ has  the smoothing effect from the Korteweg tensor.  Therefore even if we consider the conservation form \eqref{CNSK} no derivative loss occurs in the energy method for nonlinear problem in difference from CNS. 

By these linear estimate and the iteration argument in time weighted spaces, we show existence of time global solutions for small data and the decay rate of the solutions simultaneously as in \cite{Matsumura-Nishida}.

This paper is organized as follows. In section 2 notations and lemmas are described which shall be used in this paper. 
In section 3, the main result is stated. In section 4, the proof of existence of global solutions and the decay rate is stated.

%---------------------------------------------
%     Section 2.   
%------------------------------------------------------------
@

\section{Preliminaries}
In this section notations which will be used throughout this paper are introduced. Furthermore, some lemmas which will be useful in the proof of the main result are introduced.\\

The norm on $X$ is denoted by $\|\cdot\|_{X}$ for a given Banach space $X$.

Let $1\leqq p\leqq \infty.$ $L^p$ stands for the usual $L^p$ space on $\mathbb{R}^n$. 
Let $k$ be a nonnegative integer. $W^{k,p}$ and $H^k$ stand for the usual $L^p$ and $L^2$ Sobolev space of order $k$ respectively.  
(As usual, $H^0$ is defined by $H^0:=L^2$.) 

$L^p$ stands for the set of all vector fields 
$w=\trans(w_1, \cdots, w_n)$ on $\mathbb{R}^n$ 
with $w_{j}\in L^p$ $(j=1, \cdots, n)$ 
and $\|\cdot\|_{L^p}$ stands for the norm $\|\cdot\|_{(L^p)^n}$ for simplicity,  
if we have no confusion. 
Similarly a function space $X$ stands for     
the set of all vector fields 
$w=\trans(w_1, \cdots, w_n)$ on $\mathbb{R}^n$ 
with $w_{j}\in X$ $(j=1, \cdots, n)$ and  $\|\cdot\|_{X}$ stands for the norm $\|\cdot\|_{X^n}$ if no confusion will occur. 

Let $u=\trans(\phi,m)$ with $\phi\in H^k$ and $m=\trans(m_1,\cdots, m_n)\in H^j$.  
Then the norm $\|u\|_{H^k\times H^j}$ stands for  the norm of $u$ on $H^k\times H^j$, i.e., it is defined that 
$$
\|u\|_{H^k\times H^j}:=\left(\|\phi\|_{H^k}^2+\|w\|_{H^j}^2\right)^{\frac{1}{2}}.
$$
If $j=k$, $H^k$ stands for $H^k\times (H^k)^n$ for simplicity. 
The norm $\|u\|_{H^k}$ stands for 
the norm $\|u\|_{H^k\times (H^k)^n}$, i.e., it is defined that 
$$
H^k:=H^k\times (H^k)^n, 
\ \ \ 
\|u\|_{H^k}:=\|u\|_{H^k\times (H^k)^n} 
\ \ \ (u=\trans(\phi,m)).
$$
Similarly, for $u=\trans(\phi,m)\in X\times Y$ with $m=\trans(m_1,\cdots, m_n)$ 
, 
the norm $\|u\|_{X\times Y}$ stands for
$$
\|u\|_{X\times Y}:=\left(\|\phi\|_{X}^2+\|m\|_{Y}^2\right)^{\frac{1}{2}}
\ \ \ (u=\trans(\phi,m)).
$$
When $Y=X^n$, the symbol 
$X$ stands for $X\times X^n$ for simplicity,  
and $\|u\|_X$ is defined by the norm $\|u\|_{X\times X^n}$;
$$
X:=X\times X^n, 
\ \ \ 
\|u\|_{X}:=\|u\|_{X\times X^n}
\ \ \ (u=\trans(\phi,m)).
$$

%A function space with spatial weight is defined as follows. 
%For a nonnegative integer $\ell$ and $1\leq p \leq \infty$, the function space $L^{p}_{\ell}$ denotes the weighted $L^{p}$ space defined by
%$$
%L^{p}_{\ell}:=\{u\in L^{p}; \|u\|_{L^{p}_{\ell}}:=\|(1+|x|)^{\ell}u\|_{L^{p}}<\infty\}.
%$$

The symbols $\hat{f}$ and $\mathcal{F}[f]$ stand for the Fourier transform of $f$ for the space variables $x$, that is, we define that  
\begin{eqnarray}
\hat{f}(\xi)
=\mathcal{F}[f](\xi)
:=\int_{\mathbb{R}^n}f(x)e^{-ix\cdot\xi}dx\quad (\xi\in\mathbb{R}^n).\nonumber
\end{eqnarray}
In addition,  the inverse Fourier transform of $f$ is defined by
\begin{eqnarray}
\mathcal{F}^{- 1}[f](x)
:=(2\pi)^{-n}\int_{\mathbb{R}^n}f(\xi)e^{i\xi\cdot x}d\xi\quad(x\in\mathbb{R}^n).\nonumber
\end{eqnarray}

%Let $k$ be a nonnegative integer and let 
%$r_{1}$ and $r_{\infty}$ be positive constants  satisfying $r_{1}<r_{\infty}.$ 
Let $k$ be a nonnegative integer and let 
$r_{1}$ and $r_{\infty}$ be positive constants  satisfying $r_{1}<r_{\infty}.$ The symbol $H_{(\infty)}^{k}$ stands for the set of all $u\in H^k$ satisfying $\mbox{supp }\hat{u}\subset\{|\xi|\geq r_{1}\}$, 
 and the symbol $L_{(1)}^{2}$ stands for the set of all $u\in L^2$ satisfying $\mbox{supp }\hat{u}\subset\{|\xi|\leq r_{\infty}\}$. 
 
 We define operators $P_{j}$ $(j=1,\infty)$ on $L^2$ by
 \begin{eqnarray}
P_{j}f:=\mathcal{F}^{- 1}(\hat\chi_{j}\mathcal{F}[f])\quad(f\in L^2 , j=1,\infty),\nonumber
\end{eqnarray}
where 
\begin{align*}
&\hat{\chi}_{j}(\xi)\in C^{\infty}(\mathbb{R}^{n}) 
\quad(j=1,\infty),\quad 0\leq \hat\chi_{j}\leq 1 \quad(j=1,\infty),\\
&\hat\chi_{1}(\xi):=\left\{
    \begin{array}{l}
      1\quad (|\xi|\leq r_{1}),\\
      0\quad (|\xi|\geq r_{\infty}),
    \end{array}
  \right. \\ &\hat{\chi}_{\infty}(\xi):=1-\hat\chi_{1}(\xi),\\
&0<r_{1}<r_{\infty}.
\end{align*}

% and the function space $L_{(1)}^{2}$ stands for the set of all $u\in L^2$ satisfying $\mbox{supp }\hat{u}\subset\{|\xi|\leq r_{\infty}\}$. 
%It follows from Lemma \ref{lemP_1} {\rm (ii)} below that $H^k \cap L_{(1)}^2=L_{(1)}^2$ for any nonnegative integer $k$. 

% Let  $k$ and $\ell$ be nonnegative integers. The weighted $L^2$ Sobolev space $H_{\ell}^{k}$  
% is defined by 
%\begin{eqnarray}
%H_{\ell}^{k}:=\{u\in H^k;\|u\|_{H_{\ell}^{k}}< +\infty\},\nonumber
%\end{eqnarray}
%where
%\begin{eqnarray}
%\|u\|_{H_{\ell}^{k}}&:=&\left(\sum_{j=0}^{\ell}|u|_{H^{k}_{j}}^{2}\right)^{\frac{1}{2}},\nonumber\\
%|u|_{H^{k}_{\ell}}&:=&\left(\sum_{|\alpha|\leq k} \parallel |x|^{\ell}\partial^{\alpha}_{x} u\parallel _{L^2}^{2} \right)^{\frac{1}{2}}.\nonumber
%\end{eqnarray}
%Moreover, $H_{(\infty),\ell}^{k}$ denotes the weighted $L^2$ Sobolev space for the high frequency part defined by   
%H_{(\infty),\ell}^{k}:=\{u\in H_{(\infty)}^k;\|u\|_{H_{\ell}^{k}}< +\infty\}.\nonumber
%
%the weighted $L^2$ space for the low frequency part defined by 
%L^2_{(1),\ell}:=\{f\in L^2_{\ell};f\in L^2_{(1)}\}.
%
%For $-\infty\leq a<b\leq \infty$, the symbol 
%in $X$. 
%Similarly, $L^p(a,b;X)$ and $H^k(a,b;X)$ denote the $L^p$-Bochner space on $(a,b)$ 
% 
Let $0< T < +\infty$. A function space $X^s(0,T)$ stands for 
\begin{eqnarray*}
\lefteqn{X^s(0,T) =\{u=\trans(\phi,m); u=P_1 u +P_\infty u,}\\
&&P_1 u \in C([0,T]; L^2_{(1)}), \ \ P_\infty u=\trans(\phi_\infty, m_{\infty}), \\
&&\phi_\infty \in C([0,T]; H^{s+1}_{(\infty)})\cap L^2(0,T; H^{s+2}_{(\infty)}), \ \ m_\infty \in C([0,T]; H^{s}_{(\infty)})\cap  L^2(0,T; H^{s+1}_{(\infty)})\}
\end{eqnarray*}
and the norm is defined by 
\begin{eqnarray*}
\|P_1 u, P_\infty u\|_{X^s(0, T)}&=&\|P_1u\|_{C([0,T]; L^2_{(1)})}+\|P_\infty \phi\|_{C([0,T]; H^{s+1}_{(\infty)})\cap L^2(0,T; H^{s+2}_{(\infty)})}\\
&&\qquad +\|P_\infty m\|_{C([0,T]; H^{s}_{(\infty)})\cap  L^2(0,T; H^{s+1}_{(\infty)})}.
\end{eqnarray*}
Similarly, $X^s(0,\infty)$ stands for 
\begin{eqnarray*}
\lefteqn{X^s(0,\infty) =\{u=\trans(\phi,m); u=P_1 u +P_\infty u,}\\
&&P_1 u \in C([0,\infty); L^2_{(1)}), \ \ P_\infty u=\trans(\phi_\infty, m_{\infty}), \\
&&\phi_\infty \in C([0,\infty); H^{s+1}_{(\infty)})\cap L^2(0,\infty; H^{s+2}_{(\infty)}), \ \ m_\infty \in C([0,\infty); H^{s}_{(\infty)})\cap  L^2(0,\infty; H^{s+1}_{(\infty)})\}.
\end{eqnarray*}

For operators $L_1$ and $L_2$, $[L_1,L_2]$ stands for the commutator of $L_1$ and $L_2$, i.e.,
\begin{eqnarray}
[L_1,L_2]f:=L_1(L_2 f)-L_2(L_1 f).\nonumber
\end{eqnarray}

For a nonnegative number $s$, $[s]$ stands for the integer part of $s$. 

The symbol $``\ast"$ stands for the convolution on the space variable $x$. 

\vspace{2ex}

Some lemmas are stated which will be used in the proof of the main result.  

\vspace{2ex}

The following lemma is the well-known Sobolev type inequality.

\begin{lem}\label{lem2.1.} Let $n \geq 3$ and $s$ be an integer satisfying $s \geq [n/2]+1.$ Then there holds the inequality 
\begin{eqnarray}
\|f\|_{L^{\infty}}\leq C\| \nabla f\|_{H^{s-1}}\nonumber
\end{eqnarray}
for $f\in H^{s}.$
\end{lem}

%See, e.g., \cite{Galdi} for the proof of Lemma $\ref{lem Hardy anti}$.

\vspace{2ex}

The following inequalities are stated which are concerned with nonlinear estimate.

\vspace{2ex}

\begin{lem}\label{lem2.2.} 
Let $s$ be an integer satisfying $s\geq [n/2]+1$.  
Let $s_{j}$ and $\mu_{(j)}$ ($j=1,\cdots,\ell$) be nonnegative integers 
and multiindices 
satisfying $0\leq |\mu_{(j)}|\leq s_{j}\leq s+|\mu_{(j)}|$, 
$\mu=\mu_{(1)}+\cdots +\mu_{(\ell)}$, 
$s=s_{1}+\cdots+s_{\ell}\leq(\ell-1)s+|\mu|$, respectively. 
Then there holds
\begin{eqnarray}
\parallel \partial^{\mu_{(1)}}_{x}f_{1}\cdots \partial^{\mu_{(\ell)}}_{x}f_{\ell}\parallel_{L^{2}}
\leq C\prod_{1\leq j\leq \ell}\parallel f_{j}\parallel_{H^{s_{j}}}\quad  (f_j \in H^{s_j}).\nonumber
\end{eqnarray}
\end{lem}

See, e.g., \cite{Kagei-Kobayashi} for the proof of Lemma $\ref{lem2.2.}$.

\vspace{2ex}
\begin{lem}\label{lem2.3.} 
Let $s$ be an integer satisfying $s\geq [n/2]+1$. 
Suppose that $F$ is a smooth function on $I$, 
where $I$ is a compact interval of $\mathbb{R}$.  
Then for a multi-index $\alpha$ with $1\leq |\alpha|\leq s$, 
there hold the estimates 
$$
\|[\partial^{\alpha}_{x},F(f_{1})]f_{2}\|_{L^2}
\leq 
C\|F\|_{C^{|\alpha|}(I)}
\left\{1+\|\nabla f_{1}\|_{s-1}^{|\alpha|-1}\right\}
\|\nabla f_{1}\|_{H^{s-1}}\|f_2 \|_{H^{|\alpha|}},
$$ 
for $f_1\in H^s$ with $f_1(x)\in I$ for all $x\in \mathbb{R}^n$ 
and $f_2\in H^{|\alpha|}$; and 
$$
\|[\partial^{\alpha}_{x},F(f_{1})]f_{2}\|_{L^2}
\leq 
C\|F\|_{C^{|\alpha|}(I)}
\left\{1+\|\nabla f_{1}\|_{s-1}^{|\alpha|-1}\right\}
\|\nabla f_{1}\|_{H^{s}}\|f_2 \|_{H^{|\alpha|-1}},
$$
for $f_1\in H^{s+1}$ with $f_1(x)\in I$ for all $x\in \mathbb{R}^n$ 
and $f_2\in H^{|\alpha|-1}$.  
\end{lem}

See, e.g., \cite{Kagei-KawashimaCMP} for the proof of Lemma $\ref{lem2.3.}$.

\vspace{2ex}

Concerning the projections $P_1$ and $P_{\infty}$, we know the following properties. 

\vspace{2ex}

\begin{lem}\label{lemP_1}{\rm \cite[Lemma 4.2]{Okita}}
Let $k$ be a nonnegative integer. 
Then $P_{1}$ is a bounded linear operator from $L^2$ to $H^{k}$. 
In fact, it holds that 
\begin{eqnarray}
\|\nabla^{k}P_{1}f\|_{L^2}\leq C\|f\|_{L^2}\qquad (f\in L^{2}).\nonumber
\end{eqnarray}
As a result, for any $2\leq p\leq \infty$, $P_1$ is bounded from $L^2$ to $L^p$. 
\end{lem}
\vspace{2ex}

\begin{lem}\label{lemPinfty}{\rm \cite[Lemma 4.2]{Okita}, \cite[Lemma 4.4]{Kagei-Tsuda}}
{\rm (i)} 
Let $k$ be a nonnegative integer. 
Then $P_\infty$ is a bounded linear operator on $H^k$. 

\vspace{1ex}
{\rm (ii)} 
There hold the inequalities 
\begin{eqnarray*}
\|P_\infty f\|_{L^2}
& \leq &
C\|\nabla f\|_{L^2} 
\ \ (f\in H^1),
\\[2ex]
\|F_{(\infty)}\|_{L^2}
& \leq & 
C\|\nabla F_{(\infty)}\|_{L^2} 
\ \ (F_{(\infty)}\in H^{1}_{(\infty)}).
\end{eqnarray*}
\end{lem}

%------------------------------------------------------------
%     Section 2.   

%------------------------------------------------------------

\section{Main results}
In this section, a main result is stated for $(\ref{CNSK})$. We reformulate \eqref{CNSK} as follows.  
Hereafter we assume that $\rho_*=1$ without loss of generality. We set 
$\phi=\rho-1$. Substituting $\phi$ into \eqref{CNSK},
the following system is obtained; 
\begin{eqnarray}
\left\{
\begin{array}{lll}
\partial_{t}\phi +\div  m=0,\\
\partial_{t}m-\nu\Delta m-\tilde{\nu}\nabla\div m-\kappa \nabla \Delta\phi=F(u),\\
\phi|_{t=0}=\phi_0, \ \ m|_{t=0}=m_0,
\label{cnsk-nolinear}
\end{array}
\right.
\end{eqnarray}
where $u=\trans(\phi,m)$, $\nu={\mu}$, $\tilde{\nu}={\mu+\mu'}$, $\phi_0={\rho_0-1}$, 
\begin{eqnarray}
F(u)&=& -\Big\{\div((1+P_{(1)}(\phi)\phi) m \otimes m)+ \nabla (P_{(2)}(\phi)\phi^2)\\
&&\quad -\nu\Delta(P_{(1)}(\phi)\phi m)-\tilde{\nu}\nabla\div(P_{(1)}(\phi)\phi m)-\div \Phi(\phi)\Big\},\label{hisenkeikou}\\
P_{(1)}(\phi)
&=&
\int_{0}^{1}f'(1+\tau\phi)d\tau, \ \ f(\tau)=\frac{1}{\tau} \ \ (\tau \in \mathbb{R}),\nonumber\\
P_{(2)}(\phi)
&=&
\int_{0}^{1}(1-\tau)P''\Big(1+\tau\phi\Big)d\tau,\nonumber\\
\Phi(\phi)
&=&\kappa\Big\{\phi\Delta \phi I_n +(\nabla \phi)\cdot(\nabla \phi)I_n -\frac{|\nabla \phi |^2}{2}I_n-\nabla \phi \otimes \nabla \phi \Big\}.\nonumber
\end{eqnarray}
Note that \eqref{cnsk-nolinear} is not symmetric in contrast to the usual compressible Navier-Stokes system as in \cite{Matsumura-Nishida}. Therefore general theory for symmetric hyperbolic conservation law as in \cite{Kawashima1,Kawashima1} can not be applied. 

$(\ref{CNSK})$ is  linearized as follows. 
\begin{eqnarray}
\left\{
\begin{array}{lll}
\partial_{t}\phi +\div  m=0,\\
\partial_{t}m-\nu\Delta m-\tilde{\nu}\nabla\div m-\kappa \nabla \Delta\phi=0,\\
\phi|_{t=0}=\phi_0, \ \ m|_{t=0}=m_0.
\label{cnsk-linear}
\end{array}
\right.
\end{eqnarray}
By taking the Fourier transform of \eqref{cnsk-linear} with respect to the space variable $x$, we obtain the following ordinary differential equation with a parameter $\xi$. 
\begin{eqnarray}
\left\{
\begin{array}{lll}
\partial_{t}\hat{\phi}(t,\xi) +i\xi\cdot\hat{m}(t,\xi)=0,\\
\partial_{t}\hat{m}(t,\xi)+\nu|\xi|^2 \hat{m}(t,\xi)+\tilde{\nu}\xi(\xi\cdot \hat{m}(t,\xi))
+i\xi \kappa |\xi|^2\hat{\phi}(t,\xi)=0,\\
\hat{\phi}(0,\xi)=\hat{\phi}_0, \ \ \hat{m}(0,\xi)=\hat{m}_0.
\label{cnsk-f}
\end{array}
\right.
\end{eqnarray}
Therefore, the solutions of \eqref{cnsk-linear} are given by the following formulas. We define that 
$A=\displaystyle\frac{\nu+\tilde{\nu}}{2}$, $K=\displaystyle\frac{2\sqrt{\kappa}}{\nu+\tilde{\nu}}$. 
If $|\xi|\neq 0$ the Fourier transforms of $\phi$ and $m$ are given by 
\begin{eqnarray}
\hat{\phi}&=&\displaystyle\frac{\lambda_+(\xi) e^{\lambda_-(\xi) t}-\lambda_-(\xi) e^{\lambda_+(\xi) t}}
{\lambda_+(\xi) -\lambda_-(\xi)}\hat{\phi}_0-i\gamma \displaystyle\frac{ e^{\lambda_+(\xi) t}-e^{\lambda_-(\xi) t}}
{\lambda_+(\xi) -\lambda_-(\xi)}\xi\cdot \hat{m}_0,\nonumber\\
\hat{m}&=&e^{-\nu|\xi|^2 t}\hat{m}_0-i\xi\kappa|\xi|^2\left(\displaystyle\frac{ e^{\lambda_+(\xi) t}-e^{\lambda_-(\xi) t}}
{\lambda_+(\xi) -\lambda_-(\xi)}\right)\hat{\phi}_0\nonumber\\
&&\quad +\left(\displaystyle\frac{\lambda_+(\xi) e^{\lambda_+(\xi) t}-\lambda_-(\xi) e^{\lambda_-(\xi) t}}{\lambda_+(\xi) -\lambda_-(\xi)}-e^{-\nu|\xi|^2t}\right)\frac{\xi(\xi\cdot\hat{m}_0)}{|\xi|^2},\label{solution-linear}
\end{eqnarray}
where 
\begin{eqnarray}
\lambda_{\pm}(\xi)=-A(|\xi|^2\pm|\xi|^2\sqrt{1-K^2})\label{egenvalue}
\end{eqnarray}
denote roots of the characteristic equation of \eqref{cnsk-f}.

\vspace{2ex}

We obtain global $L^2$ solutions to $(\ref{CNSK})$ for small data with some regularity assumption of $m$ and decay rate of the solutions. In fact, the main result is stated as follows. 

\vspace{2ex}

\begin{thm}\label{global sol}

Let $u_0=\trans(\phi_0, m_0)\in H^{s+1}\times H^s$ where $s$ is an integer satisfying $s \geq [n/2]+1$ and $m_0 =\del_x \tilde{m}_0$. 
We also assume that $\trans(\phi_0, \tilde{m}_0)\in L^1$. We set 
$$
E_0 =\|u_0\|_{H^{s+1}\times H^s} +\|\trans(\phi_0, \tilde{m}_0)\|_{L^1}. 
$$
Then there exists a positive constant $\delta_0 $ such that if $E_0 \leq \delta_0$ there exists a global solution $u \in X^s(0,\infty)$ 
and we obtain the following decay rate of the solution: 
$$
\|\nabla^k u\|_{L^2} \leq C(1+t)^{-n/4-k/2} \ \ (k=0,1). 
$$
In addition, the uniqueness of the solution holds 
in the class 
$$
\{u=\trans(\phi,m); u=P_1 u +P_\infty u, 
\|P_1 u, P_\infty u\|_{X^s(0, \infty)}\leq  C\delta_0\}.
$$
\end{thm}

\vspace{2ex}

%\begin{rem}{\rm 
%If we do not assume that the initial data belong to $L^1$, i.e., if we suppose that $u_0=\trans(\phi_0, m_0)\in H^{s+1}\times H^s$ with $\tilde{m}_0\in L^2$, we do not obtain global $L^2$ solutions for small data when the physical case $n=3$. In fact, we have to consider a priori estimate without time weight. 
%Then in estimate of \eqref{low-nonlinear-ketsuron}, the right hand side is estimated by $C(1+t)^{1/4}$ for $n=3$ and thus depends on $t$. Hence we can not establish a priori estimate. 
%} 
%\end{rem}

%\vspace{2ex}

\section{Existence of global solutions} 
In this section, we show existence of global solutions to $(\ref{CNSK})$ and the decay rate of the solution. 
Set  
\begin{eqnarray*}
u=\trans(\phi, m), \ \ u_0=\trans(\phi_0,m_0), \ \ A=
\begin{pmatrix}
0 & \div \\
-\kappa \nabla \Delta & -\nu \Delta -\tilde{\nu}\nabla \div 
\end{pmatrix}
.
\end{eqnarray*}
Then \eqref{cnsk-nolinear} is rewritten as follows.  
\begin{eqnarray}
\del_t u+Au=f(u),  \ \ u|_{t=0}=u_0,\label{nonlinear-U}
\end{eqnarray}
where $f(u)=\trans(0,F(u))$. Based on the Duhamel principle,  we see the following integral equations. 
\begin{eqnarray}
u(t)=S(t) u_{0}+ \int_{0}^{t}S(t,\tau) F(u(\tau)) d\tau
\end{eqnarray}
where $S=S(t)$ denotes the solution operator of the system whose definition is given by \eqref{solution-linear}, that is, 
\begin{eqnarray}
S(t)\begin{pmatrix}
\phi_0 \\ 
m_0
\end{pmatrix}
=\mathcal{F}^{-1}\Big\{\begin{pmatrix}
S_{11} & S_{12} \\
S_{21} & S_{22}
\end{pmatrix}
\begin{pmatrix}
\hat{\phi}_0 \\ 
\hat{m}_0
\end{pmatrix}
\Big\},\label{sexplisit}
\end{eqnarray}
where
 \begin{eqnarray*}
S_{11}&=&\displaystyle\frac{\lambda_+(\xi) e^{\lambda_-(\xi) t}-\lambda_-(\xi) e^{\lambda_+(\xi) t}}
{\lambda_+(\xi) -\lambda_-(\xi)}, \\ 
S_{12}&=&-i \displaystyle\frac{ e^{\lambda_+(\xi) t}-e^{\lambda_-(\xi) t}}
{\lambda_+(\xi) -\lambda_-(\xi)}\trans{\xi}, \\ 
S_{21} &=& -i\xi\kappa|\xi|^2\left(\displaystyle\frac{ e^{\lambda_+(\xi) t}-e^{\lambda_-(\xi) t}}
{\lambda_+(\xi) -\lambda_-(\xi)}\right), \\ 
S_{22} &=& e^{-\nu|\xi|^2 t}I_n +\left(\displaystyle\frac{\lambda_+(\xi) e^{\lambda_+(\xi) t}-\lambda_-(\xi) e^{\lambda_-(\xi) t}}{\lambda_+(\xi) -\lambda_-(\xi)}-e^{-\nu|\xi|^2t}\right)\frac{\xi\trans \xi}{|\xi|^2}.
\end{eqnarray*}
Let 
$$
\Gamma [u] =S(t) u_{0}+ \int_{0}^{t}S(t,\tau) F(u(\tau)) d\tau
$$
To solve \eqref{nonlinear-U}, we look for a fixed point of $\Gamma$ for a given $u_0$. 
Since 
$$
P_j \Gamma [F(u)]=\Gamma [P_j F(u)] \ \ (j=1,\infty)
$$
and 
$$
\mbox{supp }\widehat{P_1 F}(u)\subset\{|\xi|\leq r_{\infty}\},
$$
$$
\mbox{supp }\widehat{P_\infty F}(u)\subset\{|\xi|\geq r_1\}
$$
for each $t\in [0,T]$, we will investigate projections of $\Gamma$ on  $L^2_{(1)}$ and $L^2_{(\infty)}$ respectively. 

\subsection{Estimates of $\Gamma $ for the low frequency part}

In this subsection, we estimate $\Gamma $ for the low frequency part. We set operators $S_1(t)$ and ${\scr S}_1(t)$ by 
$$
S_1(t)=S(t)|_{L^2_{(1)}},  \ \ {\scr S}_1(t)F = \displaystyle \int_0^t S_1(t-\tau) F (\tau)d\tau. 
$$
We show that the solution operator $S_1(t)$ is a bounded (linear) operator on $L^2_{(1)}$ for an initial data $u_{01}=\trans (\phi_{01}, m_{01})$ with $m_{01}=\del_x \tilde m_{01}$. We also show decay estimate of $S_1(t)$. In fact, we have

\vspace{2ex}

 \begin{prop}\label{S1}
{\rm (i)} 
Let $u_{01}=\trans (\phi_{01}, m_{01})$ and $m_{01}=\del_x \tilde m_{01}$. 
For each $\trans(\phi_{01}, \tilde{m}_{01})\in L^2_{(1)}$ 
and all $T>0$, $S_{1}(t)$ satisfies 
$$
S_{1}(t)u_{01} \in C([0,T];L^2_{(1)}),
$$ 
and there holds the estimate
$$
\| S_1(\cdot)u_1\|_{C([0,T];L^2_{(1)})}\leq 
C\|\trans(\phi_{01}, \tilde{m}_{01})\|_{L^2}, 
$$ 
where $T>0$ is any given positive number and $C$ is a positive constant independent of $T$. 

\vspace{1ex}

{\rm (ii)} If $\trans(\phi_{01}, \tilde{m}_{01})\in L^1$ under the assumption of {\rm (i)},  $S_{1}(t)$ satisfies the decay estimate
$$
\|\nabla^k S_1(t)u_1\|_{L^2_{(1)}}\leq 
C(1+t)^{-n/4 -k/2}\|\trans(\phi_{01}, \tilde{m}_{01})\|_{L^1} 
$$ 
for $t>0$ and $k \in \mathbb{Z}_{+}$,  where $C$ is a positive constant independent of $t$. 

\end{prop}

\vspace{2ex}

\noindent\textbf{Proof.} Due to \eqref{sexplisit}, we see that 
\begin{eqnarray}
S_1(t)\begin{pmatrix}
\phi_{01} \\ 
m_{01}
\end{pmatrix}
=\mathcal{F}^{-1}\Big\{\begin{pmatrix}
S_{11} & S_{12}\xi_j \\
S_{21} & S_{22}\xi_j
\end{pmatrix}
\begin{pmatrix}
\hat{\phi}_{01} \\ 
\hat{\tilde{m}}_{01}
\end{pmatrix}
\Big\}
\end{eqnarray}
for some $j \in \{1,\cdots n\}$. We prove {\rm (ii)} before {\rm (i)}. Concerning $S_{12}$ part in $S_1(t)$, we set 
$$
K_{12}(t, x)\tilde{m}_0=\mathcal{F}^{-1}\Big(\displaystyle\frac{ e^{\lambda_+(\xi) t}-e^{\lambda_-(\xi) t}}
{\lambda_+(\xi) -\lambda_-(\xi)}\xi_j\trans{\xi}\hat{\tilde{m}}_{01}\Big).
$$
We define a cut-off function $\chi_0=\mathcal{F}^{-1}\hat{\chi}_0$ with $\hat{\chi}_0$ satisfying 
\begin{eqnarray}
\hat{\chi}_0\in C^{\infty}(\mathbb{R}^{n}),
 \ \ 0\leq \hat{\chi}_0\leq 1, 
 \ \ \hat{\chi}_0=1 \ \ \mbox{on} \ \ \{|\xi| \leq  r_{\infty}\},
  \ \ \supp\hat{\chi}_0 \subset \{|\xi| \leq 2 r_{\infty}\}.\label{subcutofflowpart}
\end{eqnarray}
Since $\tilde{m}_{01}\in L^2_{(1)}$, we see that $\tilde{m}_{01}=\mathcal{F}^{-1}(\hat{\chi}_0\hat{\tilde{m}}_{01})$. 
We estimate $\|K_{12}(t, \cdot)\|_{L^2}$ by \eqref{egenvalue} as follows. 
\begin{eqnarray*}
\|K_{12}(t, \cdot)\|_{L^2}^2 &\leq & C\displaystyle \int_{|\xi|\leq 2r_{\infty}}\frac{|e^{\lambda_+(\xi) t}-e^{\lambda_-(\xi) t}|^2 |\xi|^4 |\hat\chi_{0}(\xi)|^2}
{|\lambda_+(\xi) -\lambda_-(\xi)|^2}d\xi\\
&\leq &C \displaystyle \int_{|\xi|\leq 2r_{\infty}}e^{-c|\xi|^2 t}d\xi \\
& \leq & C\displaystyle \int_{0}^{2r_\infty}e^{-c r^2 t} r^{n-1}dr\\
& \leq & C\Big(\frac{1}{\sqrt{t}}\Big)^n \displaystyle \int_{0}^{2\sqrt{t} r_\infty}e^{-c \eta^2}\eta^{n-1}d\eta\\
& \leq & Ct^{-n/2},
\end{eqnarray*}
where we used change of variables as $r\sqrt{t}=\eta$. Hence we get that
$$
\|K_{12}(t, \cdot)\|_{L^2} \leq  t^{-n/4}. 
$$
Similarly, we see that 
$$
\|\nabla ^k K_{12}(t, \cdot)\|_{L^2} \leq  (1+t)^{-n/4-k/2}. 
$$ 
Therefore, it holds by the Young inequality that 
$$
\|\mathcal{F}^{-1}(S_{12}\xi_j \hat{m}_0)\|_{L^2} \leq (1+t)^{-n/4-k/2}\|\tilde{m}_{01}\|_{L^1}. 
$$
Since another part of $S_1(t)$ can be estimated similarly, we have {\rm (ii)}. The estimate of {\rm (i)} also can be estimated similarly to  {\rm (ii)} and we omit the proof. Note that by definition of $S_1(t)$ and the Lebesgue convergence theorem, we obtain that $S_{1}(t)u_{01} \in C([0,T]; L^2_{(1)})$. This completes the proof.$\hfill\square$

\vspace{2ex}

We set 
$$
\Gamma_1 [F] =S_1(t) u_{01}+ {\scr S}_1(t)F
$$
for $F=\del_x f$. 
By direct application of Proposition \ref{S1}, we can derive the estimate of $\Gamma_1$.

 \vspace{2ex}
 
\begin{prop}\label{estimate-Gamma1}
{\rm (i)} 
Let $u_{01}=\trans (\phi_{01}, m_{01})$ and $m_{01}=\del_x \tilde m_{01}$. 
For each $\trans(\phi_{01}, \tilde{m}_{01})\in L^2_{(1)}$, $F=\del_x f$ with $f \in L^2(0,T; L^2_{(1)})$ 
and all $T>0$, $\Gamma_1$ satisfies 
$$
\Gamma_1 [F] \in C([0,T];L^2_{(1)}),
$$ 
and 
$$
\| \Gamma_1 [F](t)\|_{L^2}\leq  C\|\trans(\phi_{01}, \tilde{m}_{01})\|_{L^2}+C\displaystyle\int_0^t \|f(\tau)\|_{L^2}d\tau,
$$
where $C$ is a positive constant independent of $t$. 

\vspace{1ex}

{\rm (ii)} If in addition, $\trans(\phi_{01}, \tilde{m}_{01})\in L^1$ and $f \in C([0,T]; L^1)$ then $\Gamma_1$ is estimated by 
$$
\|\nabla^k \Gamma_1[F](t)\|_{L^2}\leq 
C(1+t)^{-n/4 -k/2}\|\trans(\phi_{01}, \tilde{m}_{01})\|_{L^1}+ C\displaystyle\int_0^t (1+t-\tau)^{-n/4 -k/2}\|f(\tau)\|_{L^1}d\tau
$$ 
for $t>0$ and $k \in \mathbb{Z}_{+}$,  where $C$ is a positive constant independent of $t$. 

\end{prop}

\subsection{Estimates of $\Gamma $ for the high frequency part}

In this subsection, we estimate $\Gamma $ for the high frequency part. Operators $S_\infty(t)$ and ${\scr S}_\infty(t)$ are defined by 
$$
S_\infty(t)=S(t)|_{L^2_{(\infty)}},  \ \ {\scr S}_\infty(t)F = \displaystyle \int_0^t S_\infty(t-\tau) F(\tau)d\tau. 
$$

We first show that $\{S_\infty(t)\}_{t \geq 0}$ is a $C^0$ semi-group on $H_{(\infty)}^s\times H^{s-1}_{(\infty)}$. 

 \vspace{2ex}
 
\begin{prop}\label{se-group-high} 
Let $n \geq 3$ and $s$ be an integer satisfying $s\geq [n/2]+1$. $\{S_\infty(t)\}_{t\geq 0}$ is a $C^0$  semi-group on $H_{(\infty)}^s\times H^{s-1}_{(\infty)}$ and satisfies 
$$
S_\infty(\cdot)u_{0\infty} \in C([0,T];H_{(\infty)}^s\times H^{s-1}_{(\infty)})
$$
for $u_{0\infty} \in C([0,T];H_{(\infty)}^s\times H^{s-1}_{(\infty)})$ and 
\begin{eqnarray*}
\|S_\infty(t)u_{0\infty}\|_{H_{(\infty)}^s\times H^{s-1}_{(\infty)}}\leq \|u_0\|_{H_{(\infty)}^s\times H^{s-1}_{(\infty)}}\label{sinfty}
\end{eqnarray*}
for $t \in [0,T]$, where $T>0$ is any positive number and $C$ is independent of $t$. 
\end{prop}

\vspace{2ex}

\noindent\textbf{Proof.} 
Let $F_{\infty}=\trans(F^1_{\infty},F^2_{\infty})\in H_{(\infty)}^s\times H^{s-1}_{(\infty)}$. We consider the following resolvent problem
\begin{eqnarray}
\lambda u_\infty+{A} u_\infty= F_{\infty}\label{eq:rezolent-pro}
\end{eqnarray}
for $u_\infty=\trans(\phi_\infty,m_\infty)$, where $\lambda \in \mathbb{C}$ is a parameter. 
Taking the Fourier transform of $\eqref{eq:rezolent-pro}$, we obtain 
\begin{eqnarray}
\lambda \hat u_\infty+\hat{A}_{\xi}\hat u_\infty=\hat F_{\infty},\label{eq:linear-lambda}
\end{eqnarray}
where 
\begin{eqnarray*}
\hat{A}_{\xi}=
\begin{pmatrix}
0 & i\trans{\xi}\\ 
i\kappa |\xi|\xi & \nu|\xi|^2 I_n+ \tilde{\nu}\xi \trans{\xi}
\end{pmatrix}
.
\end{eqnarray*}
Then, one can see by a similar manner  to the proof of Proposition \ref{prop-energy-high-Fourier} below that
\begin{eqnarray}
\lefteqn{\mbox{{\rm Re}}\lambda\left\{\sum_{|\alpha|=0}^{s}\Big(\kappa_1|(i\xi)^{\alpha}\hat{m}_\infty|^2+\kappa_1\kappa|(i\xi)^{\alpha}(i\xi)\hat \phi_\infty|^2+(i\xi)^{\alpha}\hat m_\infty \cdot \overline{(i\xi)^{\alpha}(i\xi)\hat {\phi}_\infty}\Big)\right\}
}\nonumber\\
& +& d_1\Big(\sum_{|\alpha|=0}^{s}|(i\xi)^{\alpha}(i\xi)\hat m_\infty|^2+\sum_{|\alpha|=0}^{s}|(i\xi)^{\alpha}|\xi|^2\hat \phi_\infty|^2\Big)\nonumber\\
&\quad& \leq C\Big\{\sum_{|\alpha|=0}^{s}|(i\xi)^{\alpha}\hat{F}^1_{\infty}|^2+\sum_{|\alpha|=0}^{s-1}(|(i\xi)^{\alpha}\hat{F}^2_{\infty}|^2\Big\}\label{resolbent2}
\end{eqnarray}
for $\xi \in \mathbb{R}^n$, where $\kappa_1$ and $d_1$ are the same  constants in \eqref{ene-claim}.  
Hence, if $\mbox{ Re} \lambda>0$, $(\lambda + \hat{A}_{\xi})^{-1}$ exists for each $\xi \in \mathbb{R}^n$ and $\hat{u}_\infty$ is represented by $\hat{u}_\infty=(\lambda + \hat{A}_{\xi})^{-1}\hat F_{\infty}$. 
We define the norm $|||\cdot|||_{s}$ on $H^s_{(\infty)}\times H^{s-1}_{(\infty)}$ by 
$$
|||u_\infty|||_{s}=\Big(\sum_{|\alpha|=0}^{s}\kappa\|\del_{x}^{\alpha} \phi_\infty\|_{L^2_{(\infty)}}^2+\sum_{|\alpha|=0}^{s-1}\|\del_{x}^{\alpha}m_\infty\|_{L^2_{(\infty)}}^2\Big)^{\frac{1}{2}}.
$$
It follows from \eqref{resolbent2} and definition of $L^2_{(\infty)}$ that 
$$
\mbox{{\rm Re}}\lambda |||u_\infty|||_{s}\leq C|||F_{\infty}|||_{s}
$$
and if $\mbox{ Re} \lambda >0$, it enjoys that 
$$
\|u_\infty\|_{H^{s}\times H^{s-1}}\leq \frac{C}{\mbox{ Re} \lambda}|||F_{\infty}|||_{s}.
$$
Hence
$$
\{\lambda;\mbox{ Re} \lambda >0\} \subset \rho (-A),
$$
where $\rho (-A)$ denotes the  resolvent set of $-A$ and it holds that 
$$
|||(\lambda +A)^{-1}F_{\infty}|||_{s}\leq \frac{C}{\mbox{ Re} \lambda}|||F_{\infty}|||_{s}.
$$
This together with the Hille-Yoshida theorem imply that $S_\infty(t)=e^{-tA}$ is a $C^0$ semigroup on $H_{(\infty)}^s\times H^{s-1}_{(\infty)}$, and we obtain \eqref{sinfty}. This completes the proof. 
$\hfill\square$

\vspace{2ex} 
Set 
$$
{\scr S}_{\infty}(t)F_{\infty}=\displaystyle\int_0^t S_\infty(t-\tau)F_{\infty}(\tau)d\tau. 
$$
By the Duhamel principle ${\scr S}_{\infty}$ is a solution operator for the linearized problem
\begin{eqnarray}
\del_t u_\infty+Au_\infty=F_{\infty},  \ \ u_\infty|_{t=0}=0\label{linear-high}
\end{eqnarray}
for $u_\infty=\trans(\phi_\infty,m_\infty)$. Furthermore, we have the following 

\vspace{2ex} 
\begin{prop}\label{prop-energy-high-Fourier}
$$
{\scr S}_{\infty}(\cdot)F_{\infty} \in C([0,T];H_{(\infty)}^{s+1}\times H^{s}_{(\infty)})\cap L^2(0,T;H_{(\infty)}^{s+2}\times H^{s+1}_{(\infty)})
$$
for $F_\infty \in L^2(0,T;H_{(\infty)}^{s}\times H^{s-1}_{(\infty)})$ and it holds that 
\begin{eqnarray}
&&\|{\scr S}_{\infty}(t)F_{\infty}\|_{H_{(\infty)}^{s+1}\times H^{s}_{(\infty)}}+\|{\scr S}_{\infty}(\cdot)F_{\infty}\|_{ L^2(0,T;H_{(\infty)}^{s+2}\times H^{s+1}_{(\infty)})} \nonumber\\
&&\qquad \leq C\|F_{\infty}\|_{L^2(0,T;H_{(\infty)}^{s}\times H^{s-1}_{(\infty)})}.\label{ene kara}  
\end{eqnarray}
for $t \in [0,T]$, where $T>0$ is any positive number and $C$ is independent of $t$ and $T$.  
\end{prop}

\vspace{2ex}
\noindent\textbf{Proof.} We use the energy estimate in the Fourier space. Our claim is to show 
\begin{eqnarray}
&&\sum_{|\alpha|=0}^{s}\displaystyle\frac{\del}{\del t}(\kappa_1\kappa\|(i\xi)^\alpha (i\xi)\hat{\phi}_{\infty}\|^2_{L^2_{(\infty)}}+\kappa_1\|(i\xi)^\alpha \hat{m}_{\infty}\|^2_{L^2_{(\infty)}})+\sum_{|\alpha|=0}^{s}\displaystyle\frac{\del}{\del t}((i\xi)^{\alpha}\hat{m}_\infty \cdot \overline{(i\xi)^{\alpha}(i\xi)\hat{\phi}_{\infty}})\nonumber\\
&&\quad +d_1\sum_{|\alpha|=0}^{s}(\|(i\xi)^\alpha |\xi|^2\hat{\phi}_{\infty}\|^2_{L^2_{(\infty)}}+\|(i\xi)^\alpha (i\xi)\hat{m}_{\infty}\|^2_{L^2_{(\infty)}})\nonumber\\
&&\qquad \leq C\|F_\infty\|^2_{H_{(\infty)}^{s}\times H^{s-1}_{(\infty)}}\label{ene-claim}, 
\end{eqnarray}
where $\kappa_1$ and $d_1$ are positive constants. 
Let $\xi \in \mathbb{R}^n$ and $|\xi|\geq r_1$. Taking the Fourier transform of \eqref{linear-high}, we see that 
\begin{eqnarray}
\left\{
\begin{array}{ll}
\displaystyle\frac{\del}{\del t}\hat{\phi}_{\infty} + i\xi \cdot \hat{m}_{\infty}=\hat{F}_{1\infty}, \\
\vspace{1ex}
\displaystyle\frac{\del}{\del t}\hat{m}_{\infty}+\nu|\xi|^2 \hat{m}_\infty+\tilde{\nu}\xi \xi \cdot \hat{m}_\infty +
i\kappa \xi |\xi|^2\hat{\phi}_{\infty} =\hat{F}_{2\infty}, \label{high-Fourier side}
\end{array}
\right.
\end{eqnarray}
where $\hat{F}_\infty =\trans(\hat{F}_{1\infty}, \hat{F}_{2\infty})$. 
For a multi-index $\alpha$ satisfying $|\alpha|\leq s$, 
taking the complex inner product of 
$(i\xi)^{\alpha}(\ref{high-Fourier side})_{2}$ 
with $(i\xi)^{\alpha}\hat{m}_{\infty}$ and taking the sum of $\alpha$ and from the real part for $|\xi| \geq r_1$ we have that 
\begin{eqnarray}
&&\frac{1}{2}\sum_{|\alpha|=0}^{s}\displaystyle\frac{\del}{\del t}|(i\xi)^{\alpha}\hat{m}_{\infty}|^2+\kappa\mbox{Re}\Big((i\xi)^{\alpha}(i\xi)|\xi|^2 \hat{\phi}_\infty \cdot \overline{(i\xi)^{\alpha}\hat{m}_{\infty}}\Big)\nonumber\\
&& \qquad + \frac{\nu}{2} \sum_{|\alpha|=0}^{s}|(i\xi)^{\alpha}(i\xi) \hat{m}_{\infty}|^2+\tilde{\nu} \sum_{|\alpha|=0}^{s}|(i\xi)^{\alpha}(i\xi) \trans \hat{m}_{\infty}|^2\nonumber\\
&&\qquad \qquad \leq C\sum_{|\alpha|=0}^{s-1}|(i\xi)^\alpha \hat{F}_{2\infty}|,\label{high-Fourier side2}
\end{eqnarray}
where $C$ is some positive constant. Note that due to \eqref{high-Fourier side} we obtain that 
\begin{eqnarray*}
(i\xi)|\xi|^2 \hat{\phi}_\infty \cdot \overline{\hat{m}_{\infty}} &=& -(|\xi|^2 \hat{\phi}_\infty \overline{i\xi \cdot \hat{m}_{\infty}})\\
&=&|\xi|^2 \hat{\phi}_\infty \overline{\displaystyle\frac{\del}{\del t}\hat{\phi}_\infty}-|\xi|^2 \hat{\phi}_\infty \overline{\hat{F}_{1\infty}}. 
\end{eqnarray*}
Therefore, we derive the inequality 
\begin{eqnarray}
&&\frac{1}{2}\sum_{|\alpha|=0}^{s}\displaystyle\frac{\del}{\del t}\Big\{|(i\xi)^{\alpha}\hat{m}_{\infty}|^2+
\kappa|(i\xi)^{\alpha}|\xi|\hat{\phi}_{\infty}|^2\Big\}
+ \frac{\nu}{2} \sum_{|\alpha|=0}^{s}|(i\xi)^{\alpha}(i\xi) \hat{m}_{\infty}|^2+\tilde{\nu} \sum_{|\alpha|=0}^{s}|(i\xi)^{\alpha}(i\xi) \cdot \hat{m}_{\infty}|^2\nonumber\\
&&\qquad \qquad \leq \epsilon \sum_{|\alpha|=0}^{s} |(i\xi)^{\alpha}|\xi|^2\hat{\phi}_{\infty}|^2
+C(\sum_{|\alpha|=0}^{s}|(i\xi)^\alpha \hat{F}_{1\infty}|^2+|\sum_{|\alpha|=0}^{s-1}|(i\xi)^\alpha \hat{F}_{2\infty}|^2), \label{enestandard}
\end{eqnarray}
where $\epsilon$ is a positive constant satisfying $\epsilon\leq \frac{c_2}{2\kappa_1}$, $\kappa_1$ and $c_2$ are defined below. 
On the other hand, we take the complex inner product of 
$(i\xi)^{\alpha}(\ref{high-Fourier side})_{2}$ 
with $(i\xi)^{\alpha}(i\xi)\hat{\phi}_{\infty}$ to obtain 
\begin{eqnarray}
&&\sum_{|\alpha|=0}^{s}\displaystyle\frac{\del}{\del t}((i\xi)^{\alpha}\hat{m}_\infty \cdot \overline{(i\xi)^{\alpha}(i\xi)\hat{\phi}_{\infty}})-(i\xi)^{\alpha}\hat{m}_\infty\cdot \overline{(i\xi)^{\alpha}(i\xi)\frac{\del}{\del t}\hat{\phi}_{\infty}}\nonumber\\
&& \quad +\nu ((i\xi)^{\alpha}|\xi|^2\hat{m}_\infty\cdot\overline{(i\xi)^{\alpha}(i\xi)\hat{\phi}_{\infty}} )
+\tilde{\nu} ((i\xi)^{\alpha}(\xi \trans\xi\hat{m}_\infty)\cdot\overline{(i\xi)^{\alpha}(i\xi)\hat{\phi}_{\infty}} )
+ \kappa \sum_{|\alpha|=0}^{s}|(i\xi)^{\alpha}|\xi|^2 \hat{\phi}_{\infty}|^2 \nonumber\\
&& \qquad \qquad =\sum_{|\alpha|=0}^{s}((i\xi)^{\alpha}\hat{F}_{2\infty}\cdot\overline{(i\xi)^{\alpha}(i\xi)\hat{\phi}_{\infty}}). \label{disipation-phi}
\end{eqnarray}
Since 
\begin{eqnarray}
\hat{m}_\infty\cdot \overline{(i\xi)\frac{\del}{\del t}\hat{\phi}_{\infty}}
=(-i\xi\cdot \hat{m}_\infty)\cdot\overline{-i\xi \cdot \hat{m}_\infty+\hat{F}_{1\infty}}
\end{eqnarray}
by $(\ref{high-Fourier side})_{2}$, we see from \eqref{disipation-phi} that 
\begin{eqnarray}
&&\sum_{|\alpha|=0}^{s}\displaystyle\frac{\del}{\del t}((i\xi)^{\alpha}\hat{m}_\infty \cdot \overline{(i\xi)^{\alpha}(i\xi)\hat{\phi}_{\infty}})
+c_2\sum_{|\alpha|=0}^{s}|(i\xi)^{\alpha}|\xi|^2 \hat{\phi}_{\infty}|^2 \nonumber\\
&&\quad \quad \leq c_3\sum_{|\alpha|=0}^{s}(|(i\xi)^{\alpha}(i\xi)\hat{m}_{\infty}|^2+|(i\xi)^{\alpha}i\xi \cdot \hat{m}_{\infty}|^2)\nonumber\\
&&\quad \qquad +C(\sum_{|\alpha|=0}^{s}|(i\xi)^\alpha \hat{F}_{1\infty}|^2+|\sum_{|\alpha|=0}^{s-1}|(i\xi)^\alpha \hat{F}_{2\infty}|^2).\label{disipation-phi-conclusion}
\end{eqnarray}
Let $\kappa_1$ be suitable large constant satisfying that $\frac{\kappa_1 \nu}{2} \geq 2c_3$, $\kappa_1 \tilde{\nu} \geq 2c_3$ and 
 \begin{eqnarray*}
 |\sum_{|\alpha|=0}^{s}((i\xi)^{\alpha}\hat{m}_\infty \cdot \overline{(i\xi)^{\alpha}(i\xi)\hat{\phi}_{\infty}})|
\leq \frac{1}{2}(\kappa_1\kappa\|(i\xi)^\alpha (i\xi)\hat{\phi}_{\infty}\|^2_{L^2_{(\infty)}}+\kappa_1\|(i\xi)^\alpha \hat{m}_{\infty}\|^2_{L^2_{(\infty)}}). 
\end{eqnarray*}
Considering $\kappa_1 \times (\ref{enestandard}) + \eqref{disipation-phi-conclusion}$ we get \eqref{ene-claim}. Integrating \eqref{ene-claim} on time and by the Plancherel theorem and Lemma \ref{lemPinfty} it holds that 
$$
\|\phi_\infty(t)\|_{H^{s+1}}^2+\|m_\infty(t)\|_{H^s}^2+d_1\displaystyle\int_0^t \|\nabla m_\infty\|_{H^s}^2+\|\nabla \phi_\infty\|_{H^{s+1}}^2 d\tau \leq C\|F_\infty\|_{L^2(0,T;H_{(\infty)}^{s}\times H^{s-1}_{(\infty)})}^2. 
$$
This implies \eqref{ene kara}. This completes the proof. 
$\hfill\square$

\vspace{2ex}
\begin{rem}
{\rm 
In the proof of Proposition \ref{prop-energy-high-Fourier}, we obtain the following energy estimate. 
\vspace{2ex} 

\begin{prop}\label{energyest}
Let $s$ be a nonnegative integer satisfying  $s\geq [n/2]+1$. Assume that 
\begin{eqnarray*}
u_{0\infty}=\trans(\phi_{0\infty},m_{0\infty})\in H^{s+1}_{(\infty)}\times H^s_{(\infty)},\\
F_\infty =\trans(F^1_\infty, F^2_\infty)\in L^2(0,T'; H^{s}_{(\infty)}\times H^{s-1}_{(\infty)})
\end{eqnarray*}
for all $T'>0$. Assume also that $u_{\infty}=\trans(\phi_{\infty},m_{\infty})$ satisfies
 \begin{eqnarray}
\left\{
\begin{array}{ll}
\partial_{t}u_{\infty}+Au_{\infty}=F_{\infty},\label{highparttimeeq}\\
u_{\infty}|_{t=0}=u_{0\infty}
\end{array}
\right.
\end{eqnarray}
and 
\begin{eqnarray*}
\phi_{\infty}\in C([0,T'];H^{s+1}_{(\infty)})\cap L^2(0,T';H^{s+2}_{(\infty)}),\ 
m_{\infty}\in C([0,T'];H^{s}_{(\infty)})\cap L^{2}(0,T';H^{s+1}_{(\infty)})
\end{eqnarray*}
for all $T'>0$. 
Then there exists an energy functional ${\cal E}[u_{\infty}]$ 
such that 
there holds the estimate 
\begin{eqnarray}\label{energy}
\frac{d}{dt}{\cal E}[u_{\infty}](t)
+{d}(\|\nabla \phi_{\infty}(t)\|_{H^{s+1}}^{2}+\|\nabla m_{\infty}(t)\|_{H^{s}}^{2})
 \leq
C\|F_{\infty}(t)\|_{H^{s}\times H^{s-1}}^2
\end{eqnarray} 
on $(0,T')$ for all $T'>0$. 
Here $d$ is a positive constant; 
$C$ is a positive constant independent of  $T'$;  
${\cal E}[u_{\infty}]$ is equivalent to $\|u_{\infty}\|_{H^{s+1}\times H^s}^2$, i.e, 
$$
C^{-1}\|u_{\infty}\|_{H^{s+1}\times H^s}^2
\leq {\cal E}[u_{\infty}]
\leq C\|u_{\infty}\|_{H^{s+1}\times H^s}^2; 
$$
and ${\cal E}[u_{\infty}](t)$ is absolutely continuous in $t\in [0,T']$ for all $T'>0$.
\end{prop}
}
\end{rem}

We set 
$$
\Gamma_\infty [F] =S_\infty(t) u_{0\infty}+ {\scr S}_{\infty}(t)F.  
$$ 
By direct application of Propositions \ref{se-group-high}-\ref{energyest}, we can derive the estimate of $\Gamma_\infty$.

 \vspace{2ex}
 
\begin{prop}\label{estimate-Gammainfty}
Let $u_{0\infty}=\trans (\phi_{0\infty}, m_{0\infty})$ and $s $ be an integer satisfying $s \geq [n/2]+1$.   
For each $\trans(\phi_{0\infty}, m_{0\infty})\in H_{(\infty)}^{s+1}\times H^{s}_{(\infty)}$, $F \in L^2(0,T; H_{(\infty)}^{s}\times H^{s-1}_{(\infty)})$ 
and all $T>0$, $\Gamma_\infty$ satisfies 
$$
\Gamma_\infty [F] \in C([0,T];H_{(\infty)}^{s+1}\times H^{s}_{(\infty)})\cap L^2(0,T;H_{(\infty)}^{s+2}\times H^{s+1}_{(\infty)}),
$$ 
and 
\begin{eqnarray}
\lefteqn{\| \Gamma_\infty [F]\|_{C([0,T];H_{(\infty)}^{s+1}\times H^{s}_{(\infty)})\cap L^2(0,T;H_{(\infty)}^{s+2}\times H^{s+1}_{(\infty)})}}\nonumber\\ 
&&\leq  C\|\trans(\phi_{0\infty}, m_{0\infty})\|_{H_{(\infty)}^{s+1}\times H^{s}_{(\infty)}}+C \|F\|_{L^2(0,T; H_{(\infty)}^{s}\times H^{s-1}_{(\infty)})},
\end{eqnarray}
where $C$ is a positive constant independent of $T$. Furthermore, $\Gamma_{\infty}[F]$ satisfies 
the estimate \eqref{energy}, i.e., 
\begin{eqnarray}\label{energy2}
\frac{d}{dt}{\cal E}[\Gamma_{\infty}[F]](t)
+{d}(\|\nabla \phi_{\infty}(t)\|_{H^{s+1}}^{2}+\|\nabla m_{\infty}(t)\|_{H^{s}}^{2})
 \leq
C\|F(t)\|_{H^{s}\times H^{s-1}}^2,\label{energy2}
\end{eqnarray} 
where $\Gamma_{\infty}[F]=\trans(\phi_{\infty},m_{\infty})$.  
\end{prop}

\vspace{2ex}

\subsection{Iteration argument to show the existence of global solutions}\label{global sol}

In this subsection, we show existence of global solutions to $(\ref{CNSK})$ for small data by the iteration argument. 
Recall that 
$$
\Gamma [u] =S(t) u_{0}+ \int_{0}^{t}S(t,\tau) F(u(\tau)) d\tau, 
$$
where $u_0=\trans(\phi_0,m_0)$, $m_0=\del_x \tilde{m}_0$ and $F(u)$ denotes the nonlinearity terms of $(\ref{CNSK})$.  
First, we denote the iteration scheme. We define $u^{(k)}$ $(k=1,\cdots)$ by 
\begin{eqnarray}
u^{(k)}=\Gamma [u^{(k-1)}]\label{ite-k}
\end{eqnarray}
and $u^{(0)}$ is given by $u^{(0)}=S(t) u_{0}$. 
Applying $P_j$ $(j=1,\infty)$ to \eqref{ite-k} respectively, we obtain that 
\begin{eqnarray}
u^{(k)}_j=\Gamma_j[u^{(k-1)}_1+u^{(k-1)}_\infty]\label{ite-k-low-high}
\end{eqnarray}
where $u^{(k)}_j=P_j u^{(k)}$, $u_{0j}=P_j u_{0}$ and $u^{(0)}_j=S_j(t) u_{0j}$. $(j=1,\infty)$. 

For any $0<T< + \infty$ we define a time weighted function space $Z^s_{a}(0,T)$ by 
\begin{eqnarray*}
\lefteqn{Z^s_{a}(0,T)=\{u=\trans(\phi,m); u=P_1 u +P_\infty u=u_1+u_\infty,}\\
&&P_1 u \in C([0,T]; L^2_{(1)}), \ \ P_\infty u=\trans(\phi_\infty, m_{\infty}), \\
&&\phi_\infty \in C([0,T]; H^{s+1}_{(\infty)})\cap L^2(0,T; H^{s+2}_{(\infty)}), \ \ m_\infty \in C([0,T]; H^{s}_{(\infty)})\cap L^2(0,T; H^{s+1}_{(\infty)}), \\
&&\|u_1,  u_\infty\|_{Z^s(0,T)} \leq a\}
\end{eqnarray*}
and the norm $\|u_1,  u_\infty\|_{Z^s(0,T)}$ is defined by 
\begin{eqnarray*}
\|u_1, u_{\infty}\|_{Z^s(0,T)}&=&\sup_{0\leq t \leq T}
\sum_{j=0}^{1}(1+t)^{\frac{n}{4}+\frac{j}{2}}\|\nabla^j  u_1\|_{L^2}\\
&&\quad +\sup_{0\leq t \leq T}(1+t)^{\frac{n}{4}+\frac{1}{2}}\| u_\infty\|_{H^{s+1}\times H^s}\\
&&\quad +\|\nabla u_\infty\|_{L^2(0,T; H^{s+1}_{(\infty)}\times H^s_{(\infty)})}\\ 
&&\quad + \sup_{0\leq t \leq T}(1+t)^{\frac{n}{4}+\frac{1}{2}}\Big(\displaystyle\int_0^t e^{-C_2(t-\tau)}(1+\tau)^{-\frac{n}{2}-1}\|\nabla u_\infty\|_{H^{s+1}\times H^{s}}^2d\tau\Big)^{\frac{1}{2}}, 
\end{eqnarray*}
where $a$ and $C_2$ are positive constants independent of $k$ and $T$ and are defined below respectively.  Note that the space $Z^s_{a}(0,T)$ has completeness with the norm $\|u_1, u_{\infty}\|_{Z^s(0,T)}$.  
By Theorem \ref{S1} {\rm (ii)}, for $j=0,1$ and $k\geq 1$ it holds that 
\begin{align}
&\|\nabla^j u^{(k)}(t)\|_{L^2}\leq \|\nabla^j S_1(t)u_{01}\|_{L^2} +\displaystyle
\int_{0}^{t}\|\nabla^j S_1(t,\tau)P_1 F(u^{(k-1)})\|_{L^2}d\tau,\label{est-nonlinear-low1}\\
&\|\nabla^j S_{1}(t)u_{01}\|_{L^2}
\leq  C(1+t)^{-\frac{n}{4}-\frac{j}{2}}\|\trans(\phi_{01}, \tilde{m}_{01})\|_{L^1}, \label{est-nonlinear-low2}
\end{align}
We estimate the second term of right hand side in \eqref{est-nonlinear-low1}. 
Due to the conservation form $F(u)$ can be represented by the divergence form, that is, 
$$
F(u)=\trans(0, \sum_{\ell=1}^{n}\del_{x_\ell} f_\ell(u))
$$
where $f_\ell(u)$ is suitable nonlinear terms given from \eqref{hisenkeikou}. (For example, each component of $m \otimes m$.)  
Hence we see from Theorem \ref{S1} {\rm (ii)}, the fact $n\geq 3$ and direct computation for $L^1$ norm of the nonlinearity that for $\|u_1^{(k-1)}, u_\infty^{(k-1)}\|_{Z^s(0,T)} \leq 1$
\begin{eqnarray}
\displaystyle
\int_{0}^{t}\|\nabla^j S_{1}(t,\tau)P_1 F(u^{(k-1)})\|_{L^2}d\tau &\leq& \displaystyle C\int_{0}^{t}
(1+t-\tau)^{-\frac{n}{4}-\frac{j}{2}}\|f_\ell (u^{(k-1)})\|_{L^1}d\tau \nonumber\\
&\leq &\displaystyle C\int_{0}^{t}(1+t-\tau)^{-\frac{n}{4}-\frac{j}{2}}(1+\tau)^{-\frac{2n}{4}}d\tau \|P_1 u^{(k-1)}, P_\infty u^{(k-1)}\|_{Z^s(0,T)}^2\nonumber\\
&\leq& C(1+t)^{-\frac{n}{4}-\frac{j}{2}}\|P_1 u^{(k-1)}, P_\infty u^{(k-1)}\|_{Z^s(0,T)}^2.\label{nonlinearityest-low3}
\end{eqnarray}
Owing to \eqref{est-nonlinear-low1}, \eqref{est-nonlinear-low2} and \eqref{nonlinearityest-low3} we get that for $k \geq 1$, $\|u_1^{(k-1)}, u_\infty^{(k-1)}\|_{Z^s(0,T)} \leq 1$ and $0<T<+\infty$ 
\begin{eqnarray}
\sup_{0\leq t \leq T}
\sum_{j=0}^{1}(1+t)^{\frac{n}{4}+\frac{j}{2}}\|\nabla^j u^{(k)}_1\|_{L^2} \leq C_0 E_0+C_1\|P_1 u^{(k-1)}, P_\infty u^{(k-1)}\|_{Z^s(0,T)}^2, \label{low-ite}
\end{eqnarray}
where constants $C_0$ and $C_1$ are independent of $k$ and $T$. Obviously, it holds that 
\begin{eqnarray*}
\sup_{0\leq t \leq T}
\sum_{j=0}^{1}(1+t)^{\frac{n}{4}+\frac{j}{2}}\|\nabla^j u^{(0)}_1\|_{L^2} \leq C_0 E_0. 
\end{eqnarray*}

Concerning estimate for $u_\infty$, we use the estimate \eqref{energy2}.  
Note that the following estimate which is related to estimate of the nonlinearity $P_\infty F(u)$ is obtained by 
direct computations based on Lemmas \ref{lem2.1.}-\ref{lemPinfty}.

\vspace{2ex} 
\begin{lem}\label{nonlinearest-high}
It holds that for $t\in [0,T]$ and $\|u_1, u_\infty\|_{Z^s(0,T)} \leq 1$ 
$$
\|P_{\infty}F(u)\|_{H^{s}\times H^{s-1}}\leq C(1+t)^{-\frac{n}{2}-1}\|u_1, u_\infty\|_{Z^s(0,T)}^2
+C(1+t)^{-\frac{n}{4}-\frac{1}{2}}\|u_1, u_\infty\|_{Z^s(0,T)}\|\nabla u_\infty\|_{H^{s+1}\times H^s}.
$$
\end{lem}

\vspace{2ex} 

\vspace{2ex}
\noindent\textbf{Proof.} We estimate $P_\infty (\phi \nabla \Delta \phi)$ which is one of the nonlinear terms. For $|\alpha|\leq s-1$ we see from Lemmas 2.1, 2.3, and 2.5 that  
\begin{eqnarray*}
\|\del_x^{\alpha }P_\infty (\phi \nabla \Delta \phi)\|_{L^2} &\leq&  \|\phi\del_x^{\alpha }\nabla \Delta \phi\|_{L^2} + \|[\del_x^{\alpha }, \phi]\nabla \Delta \phi\|_{L^2}\\
&\leq& \|\phi\|_{L^\infty}\|\del_x^{\alpha }\nabla \Delta \phi\|_{L^2} + \|\nabla \phi\|_{H^s}\|\nabla \Delta \phi\|_{H^{|\alpha|-1}}\\
&\leq& \|\nabla \phi\|_{H^{s-1}}\|\nabla \Delta \phi\|_{H^{s-1}} + \|\nabla \phi\|_{H^s}\|\nabla \Delta \phi\|_{H^{s-2}}.
\end{eqnarray*}
Hence it derives that 
$$
\|P_\infty (\phi \nabla \Delta \phi)\|_{H^{s-1}}\leq C(1+t)^{-\frac{n}{2}-1}\|u_1, u_\infty\|_{Z^s(0,T)}^2
+C(1+t)^{-\frac{n}{4}-\frac{1}{2}}\|u_1, u_\infty\|_{Z^s(0,T)}\|\nabla u_\infty\|_{H^{s+1}\times H^s}. 
$$
Since another nonlinear term can be estimated similarly, we get Lemma \ref{nonlinearest-high}. This completes the proof. 
$\hfill\square$

Let $D[u^{(k)}_\infty]=\|\nabla \phi^{(k)}_{\infty}(t)\|_{H^{s+1}}^{2}+\|\nabla m^{(k)}_{\infty}(t)\|_{H^{s}}^{2}$. 
By \eqref{energy2} and Lemma \ref{nonlinearest-high}, 
there exists a positive constant $C_2$ such that  for $t \in [0,T]$ and $\|u_1^{(k-1)}, u_\infty^{(k-1)}\|_{Z^s(0,T)} \leq 1$
\begin{eqnarray}
\lefteqn{{\cal E}[u^{(k)}_{\infty}](t)+d \displaystyle\int_{0}^{t}e^{-C_2(t-\tau)}D[u^{(k)}_\infty](\tau)d\tau}\nonumber\\
&&\quad \leq e^{-C_2 t}{\cal E}[u^{(k)}_{\infty}](0) \nonumber\\
&&\qquad  +C\|u^{(k-1)}_1, u^{(k-1)}_\infty\|_{Z^s(0,T)}^4 \displaystyle\int_{0}^{t}
e^{-C_2(t-\tau)}(1+\tau)^{-n-2}d\tau\nonumber\\
&&\qquad \quad +C\|u^{(k-1)}_1, u^{(k-1)}_\infty\|_{Z^s(0,T)}^2 \displaystyle\int_{0}^{t}
e^{-C_2(t-\tau)}(1+\tau)^{-\frac{n}{2}-1}D[u^{(k-1)}_\infty](\tau)d\tau\nonumber\\
&&\quad \leq e^{-C_2 t}{\cal E}[u^{(k)}_{\infty}](0)\nonumber\\
&&\qquad  +C\|u^{(k-1)}_1, u^{(k-1)}_\infty\|_{Z^s(0,T)}^4 (1+t)^{-n-2}\nonumber\\
&&\qquad \quad +C\|u^{(k-1)}_1, u^{(k-1)}_\infty\|_{Z^s(0,T)}^2 \displaystyle\int_{0}^{t}
e^{-C_2(t-\tau)}(1+\tau)^{-\frac{n}{2}-1}D[u^{(k-1)}_\infty](\tau)d\tau.\label{estenergy}
\end{eqnarray}
${\cal D}[u^{(k)}_\infty]$ and $\tilde{\cal E}[u^{(k)}_\infty]$ are defined by  
\begin{eqnarray*}
{\cal D}[u^{(k)}_\infty](t)&=&(1+t)^{\frac{n}{2}+1}\displaystyle\int_{0}^{t}
e^{-C_2(t-\tau)}(1+\tau)^{-\frac{n}{2}-1}D[u^{(k)}_\infty](\tau)d\tau, \\
\tilde{\cal E}[u^{(k)}_\infty](t) &=&\sup_{0\leq \tau \leq t}(1+\tau)^{\frac{n}{2}+1}{\cal E}[u^{(k)}_{\infty}](\tau)
\end{eqnarray*}
We see from \eqref{estenergy} that 
\begin{eqnarray}
\tilde{\cal E}[u^{(k)}_\infty](t)+d{\cal D}[u^{(k)}_\infty](t) &\leq& C(\tilde{\cal E}[u^{(k)}_\infty](0)+\|u^{(k-1)}_1, u^{(k-1)}_\infty\|_{Z^s(0,T)}^4+C\|u^{(k-1)}_1, u^{(k-1)}_\infty\|_{Z^s(0,T)}^2{\cal D}[u^{(k-1)}_\infty](t))\nonumber\\
&\leq& C_0^2E_0^2+ C_3^2\|u^{(k-1)}_1, u^{(k-1)}_\infty\|_{Z^s(0,T)}^4,\label{energy-est-keturon}
\end{eqnarray}
where the constant $C_3$ is independent of $k$ and $T$. 
Furthermore, due to Proposition \ref{estimate-Gammainfty} and Lemma \ref{nonlinearest-high} we derive that 
\begin{eqnarray}
\|u^{(k)}_\infty\|_{L^2(0,T;H^{s+2}\times H^{s+1})} \leq C_0E_0 +C_4\|u_1^{(k-1)}, u_\infty^{(k-1)}\|_{Z^s(0,T)}^2 \label{energy-est-keturon2}
\end{eqnarray}
when $\|u_1^{(k-1)}, u_\infty^{(k-1)}\|_{Z^s(0,T)} \leq 1$, where $C_4$ is a positive constant independent of $k$ and $T$.  
\vspace{2ex}

 Now we are in a position to prove the main result. Let $E_0 \leq \min\{\frac{1}{2C_0}, \frac{1}{4C_0C_1}, \frac{1}{4C_0C_3}, \frac{1}{4C_0C_4}\}$ and $a = \min\{1,\sqrt{C_0E_0/C_1}, \sqrt{C_0E_0/C_3}, \sqrt{C_0E_0/C_4}\}$. We see from \eqref{low-ite}, \eqref{energy-est-keturon}  and \eqref{energy-est-keturon2} that there holds that 
$$
\|P_1 u^{(k)}, P_{\infty}u^{(k)}\|_{Z^s (0,T)} \leq 2C_0E_0 \leq \min\{1, \sqrt{C_0E_0/C_1}, \sqrt{C_0E_0/C_3}, \sqrt{C_0E_0/C_4}\}, 
$$
for $u^{(k)} \in Z^s_a (0,T)$ and $k=0,\cdots$ inductively.  Furthermore, 
\begin{eqnarray}
\lefteqn{\|P_1 (u^{(k+1)}-u^{(k)}), P_{\infty}(u^{(k+1)}-u^{(k)})\|_{Z^s (0,T)}}\nonumber\\
&&\leq C_5 E_0\|u^{(k)}_1-u^{(k-1)}_1, u^{(k)}_\infty-u^{(k-1)}_\infty\|_{Z^s (0,T)}, \label{ketsuron1}
\end{eqnarray}
where $C_5 =\max\{2C_0C_1, 2C_0C_3, 2C_0C_4\}$. 
Therefore, when  in addition $E_0 \leq \frac{1}{2C_5}$ it holds that 
\begin{eqnarray}
\lefteqn{\|P_1 (u^{(k+1)}-u^{(k)}), P_{\infty}(u^{(k+1)}-u^{(k)})\|_{Z^s_a (0,T)}}\nonumber\\
&&\leq \frac{1}{2}\|u^{(k)}_1-u^{(k-1)}_1, u^{(k)}_\infty-u^{(k-1)}_\infty\|_{Z^s_a (0,T)}. \label{ketsuron2}
\end{eqnarray}
Since $T$ is any number satisfying $0<T<+\infty$ and the constants which appear in  \eqref{low-ite}, \eqref{energy-est-keturon}, \eqref{energy-est-keturon2}, \eqref{ketsuron1} and \eqref{ketsuron2} do not depend on $T$,  from the iteration 
there exists a unique global solution $u$ in the class 
$$
\{u=\trans(\phi,m); u=P_1 u +P_\infty u \in X^s(0,\infty), \|P_1u, P_\infty u\|_{X^s(0,\infty)} \leq a\}  
$$
and $u$ satisfies the decay estimate
$$
\|\nabla^k u\|_{L^2} \leq C(1+t)^{-n/4-k/2} \ \ (k=0,1). 
$$
This completes the proof. 
$\hfill\square$

\vspace{2ex}

\noindent {\bf Acknowledgements.} 
The first author is partly supported by Grants-in-Aid for Scientific Research with the Grant number: 16H03945. 
The second author is partly supported by Grant-in-Aid for JSPS Fellows with the Grant number: A17J047780. 

%------------------------------------------------------------
%     ŽQl•¶Œ£
%------------------------------------------------------------

\end{document}